\documentclass[english,12pt]{article}
\usepackage{amsxtra,amssymb,amsfonts,amsthm,amsmath}
\usepackage{graphicx}
\usepackage[latin9]{inputenc}
\usepackage{mathrsfs}
\usepackage[letterpaper]{geometry}
\usepackage{amsmath,amssymb}
\usepackage{setspace}
\usepackage{babel}
\usepackage{float}
\usepackage{cancel}

\makeatletter
\newtheorem{lemma}{Lemma}
\newtheorem{theorem}{Theorem}

\def\mathscr{\mathfrak}
\def\ds{\displaystyle}
\def\ov{\overline}
\def\R{\mathbb{R}}
\def\U{\mathcal{U}}
\def\V{\mathcal{V}}
\def\D{\mathcal{D}}
\def\G{\Gamma}
\def\l{\lambda}
\def\O{\Omega}
\def\d{\partial}
\def\o{\mathring}

\date{}

\makeatother

\begin{document}

\title{The Ritz method with Lagrange multipliers}

\author{Vojin Jovanovic\\
Systems, Implementation \& Integration\\
 Smith Bits, A Schlumberger Co.\\
 1310 Rankin Road\\
 Houston, TX 77032\\
e-mail: fractal97@hotmail.com 
\and Sergiy Koshkin\\
Computer and Mathematical Sciences\\
University of Houston-Downtown\\
One Main Street, \#S705\\
Houston, TX 77002\\
e-mail: koshkins@uhd.edu}

\maketitle

\newpage

\begin{abstract}

We develop a general form of the Ritz method for trial functions that do not satisfy the essential boundary conditions. The idea is to treat the latter as variational constraints and remove them using the Lagrange multipliers. In multidimensional problems in addition to the trial functions boundary weight functions also have to be selected to approximate the boundary conditions. We prove convergence of the method and discuss its limitations and implementation issues. In particular, we discuss the required regularity of the variational functional, the completeness of systems of the trial functions, and conditions for consistency of the equations for the trial solutions. The discussion is accompanied by a detailed examination of examples, both analytic and numerical, to illustrate the method.

\bigskip

\end{abstract}

\newpage

\section{Introduction}\label{s1}

In variational problems linear boundary conditions are often divided into essential (geometric) and natural (dynamic)
\cite[II.12]{Mikh}, \cite[4.4.7]{Reddy}. More generally, one calls the boundary conditions essential if they involve derivatives of order less than half of the order of the differential equation, and natural otherwise \cite[I.1.2]{Col}.
The common lore on the Ritz method is that the trial functions may violate the natural conditions, but must satisfy all the essential ones \cite[4.4.7]{Reddy}, \cite{Leip}. The reason is that the variational equations force the natural conditions on the trial solutions anyway, even if the trial functions themselves do not satisfy them. 

But what if we wish to use trial functions that violate the essential conditions as well? For instance, in problems involving parametric asymptotics the trial functions are pre-imposed with no regard for boundary conditions \cite{CHQZ,Gou}, and in initial-boundary problems with time-dependent boundary conditions the (time independent) trial functions can not satisfy them in principle. One may also wish to use such violating trial functions because they are simpler. Sure, it is easy enough to adjust them in one-dimensional examples, but it is not so easy at all in higher dimensions, especially in problems with fancy boundaries. Thus, there is abundant motivation to generalize the Ritz method to the trial functions that do not satisfy the essential conditions. However, surprisingly few authors consider such generalizations, many works on the subject are rather old, and sometimes give prescriptions that produce inconsistent systems and/or erroneous approximations \cite{Finl,Tiers}.

This is not to say that nothing has been done at all. One idea is to treat the essential boundary conditions as variational constraints and to remove them as any other constraints using the Lagrange multipliers. This idea is so natural that it appears occasionally in some applied works at least since 1946, see \cite{Bud}, where the authors explicitly cite the simplicity of the trial functions that violate the essential conditions as a reason for using 
them. Violating trial functions were even used in the boundary eigenvalue problems for vibrating plates, see \cite{Klein} and references therein. However, in all of these works Lagrange multipliers are essentially applied by analogy to finite dimensional optimization, numerical success serving as the only validation. A systematic discussion, or even a general description of the method, accounting for its requirements and limitations seems to be missing in the literature. This is all the more surprising since, as we shall see, the standard theory of the Ritz method \cite{DraMil,Vain,Zeid3} can be readily accommodated to handle such a generalization.

A very interesting paper \cite{Par} discusses one particular issue that arises when using the Lagrange multipliers in one-dimensional boundary value problems, completeness of the trial functions. We will discuss this and other issues that come up, extend the method to higher dimensions and give a general proof of its convergence. Our analysis indicates that more care is needed compared to the usual Ritz method, to wit, the variational functional has to be more regular on a larger space, the trial functions have to be complete in this larger space as well, and the multipliers can not be eliminated from the approximating systems using the usual variational formulas because of convergence issues. In the higher dimensional problems one needs to select boundary weight functions in addition to the trial functions, and balance the numbers of each to obtain well-posed approximating problems.

We try to keep our discussion more suggestive than technical, so the paper is structured somewhat unusually. We start in Section \ref{s2} by presenting some simple one-dimensional examples, and apply to them some natural looking approximating procedures recommended in the literature when the trial functions do not satisfy the essential boundary conditions. The examples are so simple that not only the exact solution, but even all trial solutions can be computed analytically. Distressingly, in many cases the approximations do not exist, do not converge, or converge to a wrong answer. The reasons turn out to be subtle, but they all can be traced to the vagueness in the concept of "approximation". In Sections \ref{s3}, \ref{s4} we flush out the false assumptions behind the failed approximations and develop the Ritz-Lagrange method that avoids them, proving along the way that it works. Unfortunately, in its original form the method only applies to one-dimensional problems, so more developing and proving is done in Section \ref{s5} for higher dimensions. Numerical applications to multidimensional problems follow in Section \ref{s6}. The paper ends with Conclusions, where we summarize our findings and discuss Galerkin type generalizations. Technical proofs are collected in the Appendix.

\section{One-dimensional examples}\label{s2}

In this section we test several reasonably looking methods on simple one-dimensional problems. Not all of them work and, as we will find out later, some of them work for the wrong reasons. These examples will highlight some subtleties and serve as a motivation for developing a general method later on. A typical procedure for numerical solution of boundary value problems is to represent the trial solution as a linear combination of finitely many trial functions and solve a finite dimensional problem that results. Most frequently the trial functions are required to satisfy the essential boundary conditions, otherwise some additional effort is needed to take care of them. 

\textbf{Problem 1.} Consider a boundary value problem for the second-order equation $u_{xx}=f$ on $[0,L]$ with essential boundary conditions on both ends of the interval $u(0)=u(L)=0$. We set for convenience $L=\pi$, and $f=1$ to make everything explicitly computable. The exact solution is easily found to be $\ds{\ov{u}=\frac12x^2-\frac{\pi}2x}$.

\textbf{Lanczos tau method.} This is perhaps the most popular method that uses trial functions not satisfying  essential conditions. The reader is warned that we use the broad understanding of the method, as in \cite{CHQZ,Gou} (review \cite{Finl} also describes it without using the name). Many authors, including Lanczos himself, understood it much more narrowly, only allowing Chebyshev polynomials as trial functions, see \cite{Ort}. The basic idea is to make the residual, the difference between the left and the right hand sides of the equation, orthogonal only to some trial functions used in the trial solution. This creates a shortfall in the number of equations compared to the number of unknown coefficients. This shortfall is used to impose the essential boundary conditions directly.

In our case the system to be solved can be written succinctly as 
\begin{equation}\label{tau}
\begin{cases}\ds{\int_0^L(u_{xx}-f)\delta u\,dx}\\
u(0)=u(L)=0,\end{cases}
\end{equation}
where selected trial functions are substituted for the weight $\delta u$. We select $\cos nx,\ n\geq0$ as our trial functions, they obviously do not satisfy the boundary conditions. Taking $N$ of them the trial solution is of the form 
$\ds{u^{(N)}=\frac{a_0}2+\sum_{n=1}^Na_n\cos nx}$ with unknown coefficients $a_i$ ($\frac12$ in front of $a_0$ is for agreement with the convention for the cosine series). Since there are two boundary equations we need to omit two trial functions when choosing weights, e.g. take $\delta u=\cos mx,\ m=0,\dots,N-2$ in 
Eq.\eqref{tau}. This gives 
\begin{equation}\label{costau}
\int_0^\pi(u_{xx}-1)\delta u\,dx=0
=\int_0^\pi\Big(-\sum_{n=1}^Nn^2a_n\cos nx-1\Big)\cos mx\,dx=0,
\end{equation}
where $m=0,\dots,N-2$. However, these equations are already inconsistent for any $N$ since $\int_0^\pi\cos nx\,dx=0$ for $n\geq1$ and the $m=0$ equation becomes $\int_0^\pi-1\,dx=0$. If we choose a different subset of trial functions as weights, say $m=2,\dots,N$, then the equations are consistent, but we find $-m^2a_m\frac{\pi}2=0$, i.e. $a_m=0$ for $m\geq2$. The boundary conditions then imply that also $a_0=a_1=0$ yielding $u^{(N)}=0$ for all $N$. One can see that any choice of weights here leads to inconsistency or to the trivial solution. The Lanczos tau method fails completely for our choice of trial functions. 

\textbf{Ritz method with boundary terms.} Let us turn to the Ritz method now. The corresponding variational functional is $J(u)=\int_0^L\frac12(u_x)^2+fu\,dx$ and the boundary value problem is equivalent to minimizing it on functions satisfying the boundary conditions. Since our boundary conditions are essential, and our trial functions do not satisfy them, the standard approach would have to be modified in some way. It turns out that not every plausible modification works. One obvious idea is to treat the essential conditions
as variational constraints and remove them using Lagrange multipliers. The Lagrange functional is 
$\mathscr{L}=J+\l_0u(0)+\l_\pi u(\pi)$, where $\l_0$, $\l_\pi$ are the Lagrange multipliers, 
and the variational equation is
\begin{align}\label{VarLag}
0=\delta\mathscr{L}&=\delta J+\l_0\delta u(0)+\l_\pi\delta u(\pi)+u(0)\delta\l_0+u(\pi)\delta\l_\pi\notag\\
&=\int_0^\pi u_x\delta u_x+f\delta u\,dx+\l_0\delta u(0)+\l_\pi\delta u(\pi)+u(0)\delta\l_0+u(\pi)\delta\l_\pi\notag\\
&=\int_0^\pi(-u_{xx}+f)\delta u\,dx+u_x\delta u\Big|_0^\pi
+\l_0\delta u(0)+\l_\pi\delta u(\pi)+u(0)\delta\l_0+u(\pi)\delta\l_\pi\\
&=\int_0^\pi(-u_{xx}+f)\delta u\,dx+\Big(\l_\pi+u_x(\pi)\Big)\delta u(\pi)+\Big(\l_0-u_x(0)\Big)\delta u(0)
+u(0)\delta\l_0+u(\pi)\delta\l_\pi\,.\notag
\end{align}
Since the boundary variations are independent of the internal ones we see that $\l_0=u_x(0)$ and $\l_\pi=-u_x(\pi)$.
This allows us to eliminate the multipliers from Eq.\eqref{VarLag} leading to:
\begin{equation}\label{Tiers}
0=\int_0^\pi(-u_{xx}+f)\delta u\,dx+u(0)\delta u_x(0)-u(\pi)\delta u_x(\pi)\,.
\end{equation}
This equation is derived as a generalization of the usual Ritz system e.g. in \cite{Tiers}. Let us test it on our example. As before, one substitutes the trial solution $u^{(N)}$ for $u$ and the trial functions for the weight 
$\delta u$. With our trial functions $\cos mx$ however $\delta u_x(0)=\delta u_x(\pi)=0$ for all $m$. But then 
Eq.\eqref{Tiers} reduces to the first equation in Eq.\eqref{tau} and produces the same system as in Eq.\eqref{costau}, only with $m=0,\dots,N$. We already saw that the $m=0$ equation can not be satisfied, making the entire system inconsistent for any $N$. This method does not work either.

\textbf{Ritz-Lagrange method.} Perhaps, instead of trying to eliminate the Lagrange multipliers from Eq.\eqref{VarLag} we should keep them as unknowns and supplement the obtained equations with the boundary conditions, just as one would do for any variational constraints. One positive trait is that the number of unknowns will always match the number of equations since each multiplier corresponds to a boundary condition. 
Instead of relying on Eq.\eqref{VarLag} let us find $\mathscr{L}\big(u^{(N)}\big)$ explicitly. Recall that $\ds{u^{(N)}=\frac{a_0}2+\sum_{n=1}^Na_n\cos nx}$, so
\begin{equation}\label{LaguN}
\mathscr{L}\Big(u^{(N)}\Big)=\frac{\pi}2a_0+\frac{\pi}4\sum_{n=1}^Nn^2a_n^2
+\l_0\left(\frac{a_0}2+\sum_{n=1}^Na_n\right)+\l_\pi\left(\frac{a_0}2+\sum_{n=1}^N(-1)^na_n\right)\,.
\end{equation}
Variational equations are $\ds{\frac{\d\mathscr{L}}{\d a_i}=0}$ and adding the boundary conditions we get the following system:
\begin{equation}\label{Lagsys}
\begin{cases}\pi+\l_0+\l_\pi=0\qquad\qquad\qquad\quad
\ds{\frac{a_0}2+\sum_{n=1}^Na_n=0}\\
\ds{\frac{\pi n^2}2a_n+\l_0+(-1)^n\l_\pi=0}\qquad
\ds{\frac{a_0}2+\sum_{n=1}^N(-1)^na_n=0}\,.
\end{cases}
\end{equation}
There are $N+3$ equations for $N+3$ unknowns $a_0,\dots,a_N,\l_0$ and $\l_\pi$. For even $n\neq0$ one immediately finds that $a_n=\frac2{n^2}$ from the first column equations. Analogously, for odd $n$ we have $a_n=-\frac2{\pi n^2}(\l_0-\l_\pi)$.
Then subtracting the second equation in the second column from the first
$$
\sum_{n\text{ odd}}2a_n=-\frac4{\pi}(\l_0-\l_\pi)\sum_{n\text{ odd}}\frac1{n^2}=0.
$$
This implies that $\l_0=\l_\pi=-\frac{\pi}2$ and all the odd coefficients vanish. Thus, we have 
$a_n=\frac{1+(-1)^n}{n^2}$ for $n\geq1$ and only $a_0$ remains undetermined. Adding the second column equations
$$
a_0+\sum_{n\text{ even}}2a_n=a_0+\sum_{n\text{ even}}\frac1{n^2}=0\,,\text{ and therefore}
$$
$$
a_0=a_0^{(N)}:=-4\sum_{k=1}^{\lfloor N/2\rfloor}\frac1{(2k)^2}=-\sum_{k=1}^{\lfloor N/2\rfloor}\frac1{k^2}
\xrightarrow[N\to\infty]{}-\frac{\pi^2}6\,,
$$
where $\lfloor\cdot\rfloor$ is the floor function returning the largest integer not exceeding its argument. 
Summarizing, we conclude that the trial solutions $u^{(N)}$ converge to the sum of the series 
$$
u^{(N)}(x)\xrightarrow[N\to\infty]{}u^{(\infty)}(x):=-\frac{\pi^2}{12}+\sum_{n=1}^\infty\frac{1+(-1)^n}{n^2}\cos nx\,.
$$
Recall that by extending a square integrable function $w(x)$ on $[0,\pi]$ to an even function on $[-\pi,\pi]$ one can expand it into a cosine series with coefficients
$\ds{a_n=\frac2\pi\int_0^\pi w(x)\cos nx\,dx}$ \cite[12.1]{Gold}. Considering that
\begin{equation}\label{Fuxx2}
x=\frac{\pi}{2}-\frac{2}{\pi}\sum_{n=1}^\infty\frac{1-(-1)^n}{n^2}\cos nx\quad\text{and}\quad
x^2=\frac{\pi^2}3+4\sum_{n=1}^\infty\frac{(-1)^n}{n^2}\cos nx\,
\end{equation}
we see that $u^{(\infty)}$ is exactly the cosine series for the exact solution $\ds{\ov{u}=\frac12x^2-\frac{\pi}2x}$. 

Finally, we have a method that works for this example at least. To make sure it was not an accident we will now apply it to some other boundary value problems.

\textbf{Problem 2.} Consider the biharmonic equation $u_{xxxx}=f$ with the boundary conditions $u(0)=u(\pi)=u_{xx}(0)=u_{xx}(\pi)=0$. For $f=1$ the exact solution is $\ds{\ov{u}=\frac{x^4}{24}-\pi\frac{x^3}{12}+\pi^3\frac{x}{24}}$. 

We apply the Ritz-Lagrange method again. The last two boundary conditions are natural and we need not worry about them. Thus, the variational problem is to minimize $J=\int_0^L\frac12(u_{xx})^2-fu\,dx$ on the space of functions satisfying the essential boundary conditions only. Now let us find the trial solutions. Taking $\ds{u^{(N)}=\frac{a_0}2+\sum_{n=1}^Na_n\cos nx}$ we have 
\begin{multline}\label{LaguN4}
\mathscr{L}\Big(u^{(N)}\Big)=J+\l_0u(0)+\l_\pi u(\pi)=\\
-\frac{\pi}2a_0+\frac{\pi}4\sum_{n=1}^Nn^4a_n^2
+\l_0\left(\frac{a_0}2+\sum_{n=1}^Na_n\right)+\l_\pi\left(\frac{a_0}2+\sum_{n=1}^N(-1)^na_n\right)
\end{multline}
since $\cos nx=(-1)^n$. The Ritz-Lagrange system is 
\begin{equation}\label{Lagsys4}
\begin{cases}-\pi+\l_0+\l_\pi=0\qquad\qquad\qquad
\ds{\frac{a_0}2+\sum_{n=1}^Na_n=0}\\
\ds{\frac{\pi n^4}2a_n+\l_0+(-1)^n\l_\pi=0}\qquad
\ds{\frac{a_0}2+\sum_{n=1}^N(-1)^na_n=0}\,.
\end{cases}
\end{equation}
It can be solved along the same lines as system Eq.\eqref{Lagsys} and we get $\l_0=\l_\pi=\frac{\pi}2$, 
$a_n=-\frac{1+(-1)^n}{n^4}$ for $n\geq1$ 
$a_0=a_0^{(N)}:=\frac14\sum_{k=1}^{\lfloor N/2\rfloor}\frac1{k^4}
\xrightarrow[N\to\infty]{}\frac{\pi^4}{360}$, where $\lfloor\cdot\rfloor$ is the floor function returning the largest integer less than or equal to its argument. Using Eq.\eqref{Fuxx2} and 
\begin{equation}\label{Fux3x4}
x^3=\frac{\pi^3}{4}+\sum_{n=1}^\infty\left(6\pi\frac{(-1)^n}{n^2}+\frac{12}{\pi}\,\frac{1-(-1)^n}{n^4}\right)\cos nx\,;
\ \ 
x^4=\frac{\pi^4}5+\sum_{n=1}^\infty\left(8\pi^2\frac{(-1)^n}{n^2}-48\frac{(-1)^n}{n^4}\right)\cos nx\,
\end{equation}
one can check that 
\begin{equation}\label{uinfty}
u^{(N)}(x)\xrightarrow[N\to\infty]{}u^{(\infty)}(x):=
\frac{\pi^4}{720}-\sum_{n=1}^\infty\frac{1+(-1)^n}{n^4}\cos nx
=\frac{x^4}{24}-\pi\frac{x^3}{12}+\pi^2\frac{x^2}{24}\,,
\end{equation}
and therefore $\ds{\ov{u}-u^{(\infty)}=\pi^3\frac{x}{24}-\pi^2\frac{x^2}{24}}$. This time the trial solutions do not converge to the exact one, and the limit difference is a linear combination of $x$ and $x^2$.

\textbf{Problem 3.} Consider the biharmonic equation $u_{xxxx}=f$ again, but now with all essential boundary conditions 
$u(0)=u(\pi)=u_{x}(0)=u_{x}(\pi)=0$. For $f=1$ the exact solution is $\ds{\ov{u}(x)=\frac{x^4}{24}-\pi\frac{x^3}{12}
+\pi^2\,\frac{x^2}{24}}$.

The Lagrange functional will now be $\mathscr{L}=J+\l_0u(0)+\l_\pi u(\pi)+\l'_0u_x(0)+\l'_\pi u(\pi)$, but if we use our trial functions $\cos nx$, $n\geq0$ the last two terms in $\mathscr{L}$ will be $0$ for any 
$u^{(N)}$. But then $\mathscr{L}$ is the same as in Eq.\eqref{LaguN4}, and therefore the Ritz-Lagrange system is the same as in Eq.\eqref{Lagsys4}. Then the trial solutions are also the same and they converge to $u^{(\infty)}$ from Eq.\eqref{uinfty}. However here, unlike in Problem 2, this is the right solution. The reader may wish to entertain herself by thinking over the last two problems. To help we will remove one red herring, as we shall see the presence of natural conditions is not an issue here.

\section{Continuity and convergence}\label{s3}

In this section we will analyze what went wrong (and right) in our examples and come up with a general scheme that provably works. The Lanczos tau method is the easiest to figure out. The idea behind the method is fairly intuitive. Since the exact solution has residual $0$ it is certainly orthogonal to all trial functions, and it is the only such function that also satisfies the boundary conditions. So as we construct approximations with residuals orthogonal to more and more trial functions, while also forcing the boundary conditions on them, it stands to reason that they should approach the exact solution in some sense. There are several assumptions that this argument relies on however. First, the system of trial functions should be complete, i.e. one should be able to approximate any function by their linear combinations. Second, we have to make sure that in the limit the residuals do become orthogonal to all trial functions. This was not the case in our second application of the Lanczos tau method, when $\cos mx$, $m=2,\dots,N$ were used as weights. Indeed, no residual ever had to be orthogonal to $1$ or $\cos x$. This is how the spurious solution $u=0$ slipped through the cracks with the residual $-1$. 

It is harder to explain why we got inconsistent systems when using $\cos mx$ with $m=0,\dots,N-2$. Let us use the benefit of hindsight and look at the cosine series for the exact solution $\ds{\ov{u}=-\frac{\pi^2}{12}+\sum_{n=1}^\infty\frac{1+(-1)^n}{n^2}\cos nx}$. Cutting it off at $N$ terms and computing the residual we get 
$u_{xx}-1=-\sum_{n=1}^N\big(1+(-1)^n\big)\sin nx$. Clearly, this series does not converge to $0$ in any apparent sense. But our argument above relied on exactly this kind of {\it continuity assumption: if trial solutions approach the exact solution so do the residuals, and this is not the case here.} Without the exact solution to approach 'trial solutions' have nothing to approximate, hence there is no reason for the Lanczos tau system to be consistent. And indeed it was not.

Looking closely at the above reasoning the reader will notice the vagueness in the meaning of "approximation". Approximation in what sense? To talk about approximating we need a way to measure distance between functions. Such measure is usually provided by Banach space norms. The norm relevant to the discussion above is the $L_2$ norm $\ds{\|u\|_{L_2}:=\left(\int_0^\pi|u|^2\delta u\,dx\right)^{\frac12}}$ with the corresponding inner product $\ds{\langle u,v\rangle:=\int_0^\pi uv\,dx}$. The Banach space of square integrable functions with the $L_2$ norm is a Hilbert space is
denoted $L_2([0,\pi])$, and all above references to convergence, completeness, continuity and orthogonality 
referred to this space. 

But the $L_2$ space is not sufficient for a discussion of variational methods. Many variational functionals come with a natural space of their own, the {\sf energy space} \cite[Ch.2]{Mikh}. For the functional of Problem 1 
the energy space is $\o{W}_2^1([0,\pi])$, the Hilbert space of functions square integrable with their first derivatives and vanishing at $0$ and $\pi$, with the norm $\ds{\|u\|_{W_2^1}:=\big(\|u\|_{L_2}^2+\|u_x\|^2_{L_2}\big)^{\frac12}}$. This {\sf norm is stronger} than the $L_2$ norm in the sense that any $W_2^1$ convergent sequence converges in $L_2$, but not conversely. On its energy space a variational functional $J$ typically has two important properties: it is 
{\sf continuous} and it is {\sf growing at infinity}, i.e. $J(u)\xrightarrow[\|u\|\to\infty]{}\infty$ \cite[III.10.2]{Vain}. For a convex functional (we will only consider those) these two properties are sufficient to prove convergence of the usual Ritz approximations in the energy norm \cite[6.2A]{DraMil}. Note that continuity and growth conditions balance each other: for a stronger norm continuity is preserved of course, but the functional may no longer grow at infinity, and vice versa for a weaker norm.

For our purposes the concept of energy space is not quite suitable because it hardwires essential boundary conditions are hardwired into it, and our trial functions do not satisfy them. Instead we start with a reflexive Banach space $\U$ (the reader will not lose much by assuming it to be Hilbert), that has nothing to do with the boundary conditions, and a convex functional on it $J:\U\to\R$. Next, we introduce the boundary operator, a linear map $\G:\U\to\R^s$ that maps functions into their boundary values. The subspace $\o{\U}:=\{u\in\U\,\big|\,\G u=0\}$ consists of functions that satisfy the boundary conditions. For Problem 1 we have $\,\U=W_2^1([0,\pi])$, $\G u=\big(u(0),u(\pi)\big)^T$ and $\o{\U}=\o{W}_2^1([0,\pi])$. The following three assumptions turn $\o{\U}$ into a generalized analog of the energy space:
\begin{enumerate}
\item $J:\U\to\R$ is convex and continuous;
\item $J(u)\xrightarrow[\|u\|\to\infty]{}\infty$, 
i.e. $J$ grows at infinity, on $\o{\U}$;
\item $\G:\U\to\R^s$ is linear and continuous.
\end{enumerate}
This setup applies to homogeneous boundary conditions only. Non-homogeneous boundary conditions can be accommodated in the usual manner, by selecting a function that satisfies them and switching to the differences with it. They solve the corresponding homogeneous problem and all convergence issues can be reduced to it, see e.g. \cite[2.1]{Lad}. 

The functional and the boundary conditions being dealt with we now turn to the trial functions. Recall that a system of elements in a Banach space is called complete if any element has their linear combination within any given distance from it. Let 
$\{\phi_i\}$ be a complete system in $\U$ and let $\U^{\,(N)}$ denote the linear span of $\phi_1,\dots,\phi_N$. The Ritz-Lagrange approach to approximating the minimizer of $J$ on $\o{\U}$ amounts to minimizing it on $\U^{\,(N)}$ subject to the boundary conditions. This is of course equivalent to minimizing it on $\o{\U}^{\,(N)}:=\U^{\,(N)}\cap\o{\U}$, whose elements are such linear combinations of the first $N$ trial functions that satisfy the boundary conditions. Although by assumption about completeness all functions in $\U$ can be approximated by linear combinations of the trial functions, it is not a priori clear that functions from $\o{\U}$ can be approximated by linear combinations that are themselves in $\o{\U}$. The next lemma proved in the Appendix assures us that this is the case.
\begin{lemma}\label{U0comb} For any complete system of elements in $\U$ there exists a system of their finite linear combinations belonging to $\o{\U}$ which is complete in $\o{\U}$.
\end{lemma}
This lemma effectively reduces the Ritz-Lagrange method to the traditional Ritz method. Indeed, if 
$\{\widetilde{\phi}_i\}$ is the complete system in $\o{\U}$ produced by Lemma \ref{U0comb} then applying the Ritz method with $\widetilde{\phi}_i$ as the trial functions instead of $\phi_i$ amounts to minimizing $J$ on $\o{\U}^{\,(N)}$. In other words, the Ritz-Lagrange method with $\{\phi_i\}$ produces the same $u^{(N)}$ (up to re-indexing) as the Ritz method with $\{\widetilde{\phi}_i\}$. This allows us to use well-known results on convergence of the Ritz method 
\cite[IV.12]{Vain}, \cite[6.2A]{DraMil} to prove convergence of its Ritz-Lagrange generalization.

One additional issue left is that one needs $J$ to be differentiable to minimize it using Lagrange multipliers. Note that we only need this differentiability on finite dimensional subspaces $\U^{\,(N)}$, that was explicitly true of $J$ in all our examples. In general, it suffices that $J$ is {\sf Gateaux differentiable} on $\U$, i.e. that its first variation is 
linear in $\delta u$ \cite[I.2.1]{Vain}, \cite[3.2]{DraMil}. The differentiability also automatically guarantees continuity. For general convex functionals only weak convergence can be expected, we use symbol $\xrightarrow[]{w}$ to denote it. However, due to Sobolev embedding theorems \cite[I.6]{Lad}, \cite[1.8]{Sob} weak convergence in $W_2^1$ for example implies convergence by norm in $L_2$. Summarizing we get the theorem below. The proof is fairly standard, but we outline it in the Appendix for the convenience of the reader.
\begin{theorem}\label{RitzLagCon} Suppose $J:\U\to\R$ is Gateaux differentiable, $\G:\U\to\R^s$ is a bounded linear operator, and $\{\phi_i\}$ is a complete system in $\U$. If $J$ is convex and grows at infinity on $\o{\U}$ then it has a minimizer $\ov{u}$ on it, as well as minimizers $u^{\,(N)}$ on all $\o{\U}^{\,(N)}$, and there exists a subsequence $N_k$ such that $u^{\,(N_k)}\xrightarrow[k\to\infty]{w}\ov{u}$. Moreover, if $J$ is strictly convex on $\o{\U}$ then $\ov{u},u^{\,(N)}$ are unique and $u^{\,(N)}\xrightarrow[k\to\infty]{w}\ov{u}$. In both cases the values of $J$ converge to its minimum on $\o{\U}$.
\end{theorem}
One can say more for quadratic functionals of the form $J(u)=\frac12B(u,u)+l(u)$, where $B$ is a symmetric bilinear form and $l$ is a linear form. These types of functionals produce linear boundary value problems \cite[22.1]{Zeid2a}. In Problem 1 we had $\ds{B(u,v)=\int_0^Lu_xv_x\,dx}$ and $\ds{l(u)=\int_0^Lfu\,dx}$, while in Problems 2, 3 they were $\ds{B(u,v)=\int_0^Lu_{xx}v_{xx}\,dx}$ and $\ds{l(u)=-\int_0^Lfu\,dx}$. Such $J$ are convex if $B$ is positive definite, and strictly convex if $B$ is {\sf strictly positive definite}, i.e. $B(u,u)\geq\varepsilon\|u\|^2$ for some $\varepsilon>0$. One can check that our functionals are strictly convex on $\o{\U}$ in both cases \cite[I.8]{Lad}. 
Simple algebra shows that for any $u,v\in\o{\U}$ 
\begin{equation}\label{Qvar}
J(u)-J(v)=\langle J'(v),u-v\rangle+\frac12B(u-v,u-v)\,,
\end{equation}
where $J'(v)=B(v,\cdot)+l$ is the derivative of $J$ at $v$. Hence, the first term on the right is the first variation of 
$J$ and it vanishes for a minimizer $v=\ov{u}$. Taking $u=u^{\,(N)}$ we see that for strictly positive definite $B$:
\begin{equation}\label{NormJ}
\|u^{\,(N)}-\ov{u}\|^2\leq\frac1\varepsilon\,B(u^{\,(N)}-\ov{u},u^{\,(N)}-\ov{u})
=\frac2\varepsilon\Big(J(u^{\,(N)})-J(\ov{u})\Big)\xrightarrow[N\to\infty]{}0.
\end{equation}
Thus, the Ritz-Lagrange solutions converge to the minimizer even by norm. More general conditions for convergence by norm are given in \cite[6.2A]{DraMil}.
 
Theorem \ref{RitzLagCon} mostly justifies the Ritz-Lagrange method used to solve Problem 1. 
The required properties of $J(u)=\int_0^L\frac12(u_x)^2+fu\,dx$ and the boundary operator $\G(u)=\big(u(0),u(L)\big)^T$ are easily verified, except for the strict convexity. That one follows from the Poincar\`e-Friedrichs inequality \cite[I.6]{Lad}. We will postpone the discussion of one other issue until the next section, and focus here on the reason why the Ritz method with boundary terms did not work. At first glance, it seems to differ from the Ritz-Lagrange method only in the manner the system for $u^{\,(N)}$ is derived. We would indeed get the same system with the same solutions if we used Eq.\eqref{VarLag} as the starting point rather than Eq.\eqref{Tiers}. The innocent step of replacing Lagrange multipliers with their values seems to have spoiled the outcome (the reader is welcome to stop reading here and think about this for a moment).

It was not so innocent after all. It helps to use the hindsight again. We have $\l_0=\ov{u}_x(0)$ and $\l_\pi=-\ov{u}_x(\pi)$ for the exact solution, but for any finite $N$ the trial solutions have $u^{(N)}_x(0)=u^{(N)}_x(\pi)=0$, while  the Lagrange multipliers are $\l_0=\l_\pi=-\frac{\pi}2\neq0$. We see that {\it the values of derivatives of the trial solutions at the ends of the interval do not converge to the values of the Lagrange multipliers, even though the two are equal for the exact solution}. Substitution of the limit values in Eq.\eqref{Tiers} relies on exactly the same kind of hidden continuity assumption that was made about the residual in the Lanczos tau method. 

How does this reconcile with convergence of the Ritz-Lagrange trial solutions to the exact solution? We do have convergence but even for strictly positive definite quadratic functionals it is not strong enough to provide such continuity. Indeed, convergence by norm in $W_2^1$ means that derivatives $u^{(N)}_x$ converge to $\ov{u}_x$ in $L_2$, not pointwise. In fact, we will always have $u^{(N)}_x(0)=u^{(N)}_x(\pi)=0$ regardless of the variational problem as long as we use cosines for trial functions. Unless it so happens that the first derivatives of the exact solution vanish at the endpoints there will never be convergence of the end point derivatives to the Lagrange multipliers. {\it Lagrange multipliers have to be kept as variables and can not be assumed to satisfy relations that hold for the exact solution}.

\section{Completeness}\label{s4}

Nothing we discussed so far explains why the Ritz-Lagrange method did not work for Problem 2. For this problem 
$\U=W_2^2([0,\pi])$, the Hilbert space of functions square integrable with their first and second derivatives, with the norm $\ds{\|u\|_{W_2^2}:=\big(\|u\|_{L_2}^2+\|u_x\|^2_{L_2}+\|u_{xx}\|^2_{L_2}\big)^{\frac12}}$. The functional is 
$J(u)=\int_0^L\frac12(u_{xx})^2-fu\,dx$ and the boundary operator is again $\G(u)=\big(u(0),u(L)\big)^T$. 
The space $\o{\U}$ consists of functions from $W_2^2([0,\pi])$ that vanish at the ends, this space is denoted 
$W_{2,0}^2([0,\pi])$ in \cite[II.6]{Lad}. Problem 3 uses the same functional but $\G(u)=\big(u(0),u(L),u_x(0),u_x(L)\big)^T$, so $\o{\U}=\o{W}_2^2([0,\pi])$, the space of $W_2^2$ functions that vanish at the ends along with their derivatives. One can check that in both cases $J$ and $\G$ satisfy all the conditions of Theorem \ref{RitzLagCon}. And yet for Problem 2 we did not obtain the exact solution in the limit. 

After inspecting the theorem closely the reader will notice one more condition that we did not address yet, because (perhaps) it seemed to be obviously satisfied. But "when you have eliminated the impossible, whatever remains, however improbable, must be the truth". That would be the completeness of trial functions, and \cite{Par} deserves credit for highlighting how far from obvious it can get. Theorem \ref{RitzLagCon} requires $W_2^2$ completeness in this case. What is obvious, or at least well-known from the standard theorems on the Fourier series, is that any function in $L_2([0,\pi])$ can be approximated by cosines (sines) in the $L_2$ norm. This is because any $L_2$ function on $[0,\pi]$ can be extended by evenness (oddness) to $[-\pi,\pi]$ while remaining in $L_2$. But that is not even enough for Problem 1, where $W_2^1$ completeness was required (this is the one issue we skipped). Fortunately, $W_2^1$ completeness of cosines reduces to the $L_2$ completeness of sines. 
\begin{lemma}\label{CosW21} The system $\cos nx$, $n\geq0$ is complete and minimal in $W_2^1([0,\pi])$.
\end{lemma}
\noindent The {\sf minimality} above means that the system becomes incomplete after deleting any function, the proof is in the Appendix. 

Not so for $W_2^2$, {\it cosines are incomplete in $W_2^2([0,\pi])$}. As observed in \cite{Par}, the second derivatives $0, -\cos x, -4\cos 2x,\dots,-n^2\cos nx,\dots$ do not include a constant, and therefore can not approximate the second derivative of $x^2$ in $L_2$. But then cosines can not approximate $x^2$ in $W_2^2$ since its norm incorporates the $L_2$ norm for second derivatives. In \cite{Par} the authors add $x^2$ to the cosine system, but... we shall see that adding 
$x^2$ is not enough. Here is a quick way to see it. In contrast to $W_2^1$, pointwise values of first derivatives are well-defined and continuous on $W_2^2$ (this follows from the Sobolev embedding theorems \cite[I.8]{Lad}, 
\cite[1.8]{Sob}). Since cosines satisfy $u_x(0)=u_x(\pi)=0$ any function that can be approximated by them must satisfy the same equalities. These are two independent conditions, so the codimension of cosines' linear span is at least two, while $x^2$ spans only one extra dimension. It follows from the next Lemma that the codimension is exactly two and $x$ also needs to be added. 
\begin{lemma}\label{CosW22} The system $x,x^2,\cos nx$, $n\geq0$ is complete and minimal in $W_2^2([0,\pi])$.
\end{lemma}
It may seem odd that we have to add $x$ on top of $x^2$, after all there is no second derivative issue with it. However, while $L_2$ approximation of the second derivatives is necessary for $W_2^2$ approximation of functions themselves, it is not sufficient. {\it We need to approximate the second derivative and the function itself using the same expression.} The second derivative of $x$ can most certainly be approximated by cosines in $L_2$, as can be $x$ itself, but not while the former approximation is the second derivative of the latter. In particular, the cosine series for $x$ diverges in $L_2$ after two differentiations, see Eq.\eqref{Fuxx2}. Note that {\it verifying completeness in a correct space can not be avoided even if one uses the usual Ritz method with trial functions satisfying the essential boundary conditions.}

What Problem 2 demonstrates is that solutions can not always be approximated in the space where the trial functions are complete. Only completeness in the norm dictated by the variational functional counts. {\it Completeness of trial functions in a norm weaker than the energy norm does not simply weaken the convergence to the exact solution, trial solutions may not converge to it at all.} Following Lemma \ref{CosW22}, we should add $x$ and $x^2$ as the trial functions and represent the trial solutions as
$\ds{u^{(N)}=b_1x+b_2x^2+\frac{a_0}2+\sum_{n=1}^Na_n\cos nx}$. The Lagrange functional becomes
\begin{multline}\label{LaguN4x}
\mathscr{L}\Big(u^{(N)}\Big)=-\frac{\pi^2}2b_1-\frac{\pi^3}3b_2-\frac{\pi}2a_0+2\pi b_2^2+\frac{\pi}4\sum_{n=1}^Nn^4a_n^2\\
+\l_0\left(\frac{a_0}2+\sum_{n=1}^Na_n\right)+\l_\pi\left(\pi b_1+\pi^2b_2+\frac{a_0}2+\sum_{n=1}^N(-1)^na_n\right)\,,
\end{multline}
and the Ritz-Lagrange system is
\begin{equation}\label{Lagsys4x}
\begin{cases}\ds{-\frac{\pi^2}2+\pi\l_\pi=0}\qquad\qquad\qquad\ 
\ds{\frac{\pi n^4}2a_n+\l_0+(-1)^n\l_\pi=0}\\
\ds{-\frac{\pi^3}3+4\pi b_2+\pi^2\l_\pi=0}\qquad\quad\!
\ds{\frac{a_0}2+\sum_{n=1}^Na_n=0}\\
-\pi+\l_0+\l_\pi=0\qquad\qquad\qquad\!\!
\ds{\pi b_1+\pi^2b_2\frac{a_0}2+\sum_{n=1}^N(-1)^na_n=0}\,.
\end{cases}
\end{equation}
The equations with $\l$-s immediately yield $\l_0=\l_\pi=\frac{\pi}2$, $b_2=-\frac{\pi^2}{24}$ and 
$\ds{a_n=-\frac{1+(-1)^n}{n^4}}$, $n\ge1$. From the second equation in the second column 
$$
a_0=a_0^{(N)}:=2\sum_{n=1}^{N}\frac{1+(-1)^n}{n^4}=\frac14\sum_{k=1}^{\lfloor N/2\rfloor}\frac1{k^4}
\xrightarrow[N\to\infty]{}\frac{\pi^4}{360}\,
$$
as before. Finally, subtracting the second equation in the second column from the third we have $b_1=-\pi b_2=\frac{\pi^3}{24}$. Thus,
\begin{gather*}
u^{(\infty)}(x)=\frac{\pi^3}{24}x-\frac{\pi^2}{24}x^2+\frac{\pi^4}{720}-\sum_{n=1}^\infty\frac{1+(-1)^n}{n^4}\cos nx
=\frac{\pi^3}{24}x-\cancel{\frac{\pi^2}{24}x^2}
+\frac{x^4}{24}-\pi\frac{x^3}{12}+\cancel{\pi^2\frac{x^2}{24}}=\ov{u}(x)
\end{gather*}
matches the exact solution as expected.

This does it for Problem 2, but now we get an unexpected puzzle in Problem 3. How come we got the correct answer while using an incomplete system of trial functions? If a system $\{\phi_i\}$ is incomplete in $\U$ the arguments leading to 
Theorem \ref{RitzLagCon} still apply to the closure of their linear span $\U_\phi$ in place of $\U$. Therefore, the Ritz-Lagrange approximations converge to the minimizer $\ov{u}_\phi$ on $\o{\U}_\phi:=\o{\U}\cap\U_\phi$ rather than $\o{\U}$. For quadratic functionals there is a simple relation between $\ov{u}$ and $\ov{u}_\phi$. 
By Eq.\eqref{Qvar} we have $\ds{J(u)=J(\ov{u})+\frac12B(u-\ov{u},u-\ov{u})}$ and minimizing $J$ on 
$\o{\U}_\phi$ is equivalent to minimizing $B(u-\ov{u},u-\ov{u})$ there. As $B$ is a quadratic form, the solution is well known to be the orthogonal projection of 
$\ov{u}$ to $\o{\U}_\phi$ with respect to the inner product defined by $B$ \cite[IV.11]{CH}. Recall from Lemma \ref{CosW22} that the only functions missing from the system of cosines in $W_2^2$ are $x$ and $x^2$. The exact solution $\ds{\ov{u}(x)=\frac{x^4}{24}-\pi\frac{x^3}{12}+\pi^2\,\frac{x^2}{24}}$ is $B$-orthogonal to both of them, despite the misleading appearance of $x^2$ in the formula, as one can check by direct integration. Thus, the missing functions are not needed to approximate $\ov{u}$ simply because $\ov{u}$ happens to be in (the $W_2^2$ closure of) the linear span of cosines. This can also be seen directly from its cosine series $\ds{\ov{u}=\frac{\pi^4}{720}-\sum_{n=1}^\infty\frac{1+(-1)^n}{n^4}\cos nx}$, which converges not only in $L_2$ but even in $W_2^2$.

\section{The Ritz-Lagrange method in higher dimensions}\label{s5}

The Ritz-Lagrange method  described in Section \ref{s2} can not be applied to multidimensional boundary value 
problems. In this section we will develop a suitable generalization and prove that it works. The main distinction is that the boundary operator $\G:\U\to\V$ no longer maps into a finite dimensional space. Indeed, in dimensions two and higher the boundary values are not arrays of numbers, but functions on the boundary forming an infinite dimensional space $\V$. The induction proof of the key Lemma \ref{U0comb} no longer works and its claim itself is false. It is easy to find complete systems of functions with no (finite) non-trivial linear combinations satisfying the boundary conditions. {\it If we are committed to using arbitrary complete systems of trial functions we must find a way to form their linear combinations that satisfy essential boundary conditions "approximately"}.

To this end we will use a complete system $\{\psi_j\}$ of linear functionals on $\V$, i.e. elements of the dual space 
$\V^*$ (as with $\U$, the reader may assume that $\V$ is a Hilbert space, in which case $\V^*=\V$). If 
$\langle \psi_j, \G u\rangle=0$ for all $j$ then $\G u=0$ and $u\in\o{\U}$, so we can think of operators 
$\G_s(u)=\big(\langle \psi_1, \G u\rangle,\dots\langle\psi_s, \G u\rangle\big)^T$ as approximations to $\G$, and the corresponding spaces $\,\o{\U}_s:=\{u\in\U\,\big|\,\G_s(u)=0\}$ as approximations to $\o{\U}$. Assuming $\G$ is 
continuous $\G_s$ will be also and we can apply Theorem \ref{RitzLagCon} with $\G_s$ in place of $\G$ for each $s$. This gives us a sequence of approximations $u_s^{(N)}$ converging to an exact minimizer $\ov{u}_s$ of $J$ on $\o{\U}_s$. The remaining question is whether we can count on $\ov{u}_s$ to approximate the overall minimizer $\ov{u}$ of $J$ on 
$\o{\U}$. Before proceeding let us describe the approximating procedure that our approach suggests. 
\pagebreak
\begin{quote}
\textbf{Multidimensional Ritz-Lagrange method.}\\
To minimize a functional $J:\U\to\R$ subject to essential boundary conditions $\G u=0$ with $\G:\U\to\V$ select {\sl internal trial functions} $\phi_1,\dots,\phi_N\in\U$ and {\sl boundary weight functions} $\psi_1,\dots,\psi_s\in\V^*$ with $N\gg s$. A {\sl Ritz-Lagrange trial solution} $\ds{u_s^{(N)}=\sum_{i=1}^Na_i\phi_i}$ is obtained by solving the system of $N+s$ equations with $N+s$ unknowns $a_1,\dots,a_N,\l_1,\dots,\l_s$ consisting of $N$ {\sl internal equations} 
$\ds{\frac{\partial\mathscr{L}}{\partial a_i}=0}$ and $s$ {\sl boundary equations} $\langle \psi_j, \G u\rangle=0$, where 
$\mathscr{L}:=J(u_s^{(N)})+\langle\l^{(s)},\G u_s^{(N)}\rangle$ is the Lagrange functional and 
$\l^{(s)}:=\l_1\psi_1+\dots+\l_s\psi_s$ is the Lagrange multiplier.
\end{quote}
For higher order equations it is more natural to use several boundary operators $\G_k$ instead of a single one, e.g. one for function values, another for normal derivatives, etc. One can always formally wrap them into a single operator with combined codomain, so no additional theoretical discussion is necessary. 

A justification of our method is based on Theorem \ref{RitzLagCon} and Theorem \ref{RitzLagBnd} below. The reader not interested in justification may skip the rest of this section and look at applications in the next one. One piece of bad news is that we can not expect $\ov{u}_s$ to converge to $\ov{u}$ in the same generality as in Theorem \ref{RitzLagCon}. The root cause is that {\it the minimizer in $\,\o{\U}$ is approximated by elements outside of $\,\o{\U}$, which is why we need $J$ to be well-behaved on the entire $\U$, not just $\,\o{\U}$,} something avoidable if $\phi_i$ do satisfy the boundary conditions. In particular, {\it the values $J(\ov{u}_s)$ are potentially smaller than $J(\ov{u})$ because they are obtained by minimizing $J$ on larger subspaces $\,\o{\U}_s\supset\o{\U}$}. As a consequence, standard properties of convex functionals, that we relied on in Theorem \ref{RitzLagCon}, no longer guarantee convergence of $J(\ov{u}_s)$ to $J(\ov{u})$ if $\ov{u}_s$ converges only weakly (in technical terms, continuous convex functionals are weakly lower semi-continuous, but usually not weakly upper semi-continuous 
\cite[III.8.5]{Vain}, \cite[41.2]{Zeid3}). To make our proof work we need to assume a stronger form of convexity. For  Gateaux differentiable functionals $J$ convexity is equivalent to monotonicity of their derivatives, i.e. 
$\langle J'(u)-J'(v),u-v\rangle\geq0$ for all $u,v$ \cite[II.5.3]{Vain}. This is a generalization of a familiar fact that convex functions have monotone derivatives. We will need a form of {\sf uniform monotonicity}, cf. \cite[V.18.6]{Vain}, namely
\begin{equation}\label{Umon}
\langle J'(u)-J'(v),u-v\rangle\geq c\,(\|u-v\|),
\end{equation}
where $c\,(t)$ is a continuous monotone increasing function with $c(t)=0$. The point is that if $c\,(\|u_s\|)\to0$ then 
$\|u_s\|\to0$ and hence $u_s\to0$ by norm. We are now ready to state the main theorem of this section.
\begin{theorem}\label{RitzLagBnd} Suppose $J:\U\to\R$ is Gateaux differentiable and $\G:\U\to\V$ is a bounded linear map.
Let $\{\psi_j\}$ be a complete system in $\V^*$ and set $\,\o{\U}_s:=\{u\in\U\,\big|\,\langle \psi_j, \G u\rangle=0,\,1\leq j\leq s\}$. If $J$ grows at infinity on some $\o{\U}_{s_0}$, and its derivative on it satisfies Eq.\eqref{Umon} 
then it has a minimizer $\ov{u}$ on $\o{\U}$, as well as minimizers $\ov{u}_s$ on all $\o{\U}_s$ with $s\geq s_0$, and $\ov{u}_{s}\xrightarrow[s\to\infty]{}\ov{u}$. The values of $J$ converge to its minimum on $\o{\U}$.
\end{theorem}
In examples it is typical that $J$ does not satisfy Eq.\eqref{Umon} on the entire space $\U$, but does satisfy on subspaces much larger than $\o{\U}$, such as $\o{\U}_{s}$. Let us discuss the case of quadratic functionals 
$J(u)=\frac12B(u,u)+l(u)$ in more detail. From the calculation in Eq.\eqref{Qvar} we know that 
$J'(v)=B(v,\cdot)+l$. Therefore,
\begin{equation*}\label{Bmon}
\langle J'(u)-J'(v),u-v\rangle=B(u,u-v)-B(v,u-v)=B(u-v,u-v)\,.
\end{equation*}
We need $B(u,u)\geq\varepsilon\|u\|^2$ with $\varepsilon>0$ to satisfy Eq.\eqref{Umon}, i.e. we need $B$ to be strictly positive definite. Take the multidimensional analog of the functional from Problem 1 as an example, $J(u)=\int_{\O}\frac12(\nabla u)^2+fu\,dx$, where $\O$ is a domain with smooth boundary. It follows from the Poincar\`e-Friedrichs inequality \cite[I.6]{Lad}, \cite{Zeid2a} that $B(u,u)=\int_{\O}\frac12(\nabla u)^2\,dx$ is strictly positive definite on 
$\o{W}_2^1(\O)$, but it most certainly is not on the entire $W_2^1(\O)$ since $B(u,u)=0$ for any $u=\text{const}$. Nevertheless, it still follows from the calculus of variations that $B$ satisfies Eq.\eqref{Umon} on any subspace complementary to the constants, see e.g. \cite[VI.1]{CH}. Similar considerations apply to other quadratic forms related to the strongly elliptic equations like the biharmonic equation. They are usually strictly positive definite on complements to finite-dimensional subspaces that they annihilate \cite[22.11]{Zeid2a}.

It is worth stressing that {\it Theorems \ref{RitzLagCon}, \ref{RitzLagBnd} do not imply that the double sequence $u_s^{(N)}$ converges to $\ov{u}$}. In fact, let $N<s$ and let the functionals $\G^*\psi_1,\dots,\G^*\psi_s$, where $\G^*:\V^*\to\U^*$ is the dual of $\G$, be linearly independent. Then the boundary equations 
$\langle \psi_j, \G u\rangle=\langle\G^*\psi_j,u\rangle=0$ alone are enough to force $u_s^{(N)}=0$ no matter how large 
$N$ and $s$ are. In practice, this means that {\it one should always take many more internal trial functions than the boundary ones}, hence $N\gg s$. This way for large $N$ the approximation $u_s^{(N)}$ will be close to $\ov{u}_s$ by Theorem \ref{RitzLagCon}, while $\ov{u}_s$ in turn will be close to $\ov{u}$ by Theorem \ref{RitzLagBnd} if $s$ itself is large enough. 

As in one-dimensional examples one will have to verify completeness of trial functions, both internal and boundary, in the appropriate space. One has to be extra careful with functionals involving higher order derivatives because the values of function and their derivatives have to be approximated simultaneously. Natural spaces to use are $W_p^k(\O)$, the spaces of functions with integrable $p$-th powers along with all of their derivatives up to order $k$. 
A generalization of the Weierstrass theorem implies that polynomials form a complete system in $W_p^k(\O)$
for any $p\geq1$, any bounded domain $\O$, and any positive integer $k$ (in fact, polynomials are even uniformly complete \cite[II.4.3]{CH}). However, polynomials may not always be convenient in a particular problem. The following Lemma 
can be useful in finding other complete systems. 
\begin{lemma}\label{MultDComp} Let $\{\phi_i\}$ and $\{\widetilde{\phi}_j\}$ be complete systems in $W_p^k(\O)$ and $W_p^k(\widetilde{\O})$ respectively, where $\O$ and $\widetilde{\O}$ are some bounded domains. Then the system 
$\{\phi_i\widetilde{\phi}_j\}$ is complete in $W_p^k(\O\times\widetilde{\O})$.
\end{lemma}
\noindent If one starts from one-dimensional systems the lemma will only produce complete systems in box-like domains $[a_1,b_1]\times\cdots\times[a_m,b_m]$. However, any system of functions complete on a domain will be complete on any of its subdomains, so for an arbitrary domain one can always use a system complete on the smallest box that contains it. Unfortunately, such product systems can not be expected to be minimal even if the original systems were. 
A more targeted choice is to take eigenfunctions of an operator on the same domain that is simpler than the one involved, but is somewhat similar to it. Various spectral theorems often ensure completeness of eigenfunctions in suitable Sobolev spaces \cite[22.11a]{Zeid2a}.

\section{Multidimensional examples}\label{s6}

In this section we illustrate the multidimensional Ritz-Lagrange method developed in Section \ref{s5} by applying it to some typical problems. Since calculations by hand quickly become intractable we performed them using a computer algebra system.

\textbf{Problem 4.}  Consider a boundary value problem for the Laplace equation $\nabla^2u=f$ in $\Omega$, where $\Omega$ is the unit disk, with the boundary condition $u=0$ on $\partial{\Omega}$ and $f=\cos(\sqrt{x^2+y^2})$. This equation describes the transverse deflection of a membrane fixed everywhere at the boundary and subjected to pressure given by $f$ \cite[IV.10.3]{CH}. 

The profile of $f$ was chosen so that the problem has an analytic solution which is not a polynomial. Specifically, one can represent the exact solution as a rapidly convergent series
\begin{equation} \label{CiSol}
\ov{u}(r)=\gamma+\cos(1)-\mathrm{Ci}(1)-\cos(r)
+\sum_{i=1}^\infty{\frac {(-1)^{i}}{2i \left( 2\,i \right) !}}\,{r}^{2\,i} 
\end{equation}
where $r=\sqrt{x^2+y^2}$, $\gamma:=\lim_{n\to\infty}\Big(\sum_{k=1}^n\frac1k-\ln n\Big)$ is the Euler-Mascheroni constant, and $\mathrm{Ci}(x):=-\int_x^\infty\frac{\cos t}{t}\,dt$ is the cosine integral.
 
To solve this problem we use the multidimensional Ritz-Lagrange method. The variational problem is to minimize the  functional $\ds{J(u)=\int_\Omega \frac{1}{2}(\nabla u)^2+fu\,dxdy}$, which gives the total potential energy of the membrane, subject to the boundary condition. In the notation of Section \ref{s5} we take $\U=W_2^1(\O)$ with $\G$ being the restriction of $u$ to the boundary $\partial{\Omega}$. Moreover, $\G$ is continuous if we take $\V=L_2(\partial{\Omega})$. Our internal trial functions are the monomials, which obviously do not satisfy 
the boundary condition, and the trial solution is $\ds{u^{(N)}=\sum_{i=1}^N \sum_{j=1}^N  c_{ij}\,x^{i-1}y^{j-1}}$. Note that $N$ of Section \ref{s5} will be $N^2$ here because of double indexing. As the boundary weight functions we choose the piecewise linear ones on uniform partitions of $\partial{\Omega}$. Unlike the usual choices for a circle, e.g. the trigonometric functions, these ones can be used on a wide variety of boundaries. Instead of using a single indexed system $\psi_k$ it is convenient to split it into the constant $g_k$ and the linear $h_k$ parts. If the boundary is partitioned into $s$ segments we have 
\begin{equation} \label{gh_functions}
g_{k}(\theta):=\begin{cases}
\ds{1,\,\frac{2\pi k}{s} \le \theta \le \frac{2\pi(k+1)}{s}}\\
0,\ \mbox{otherwise}
\end{cases}\text{and}\quad 
h_{k}(\theta):=\begin{cases}
\ds{\theta,\,\frac{2\pi k}{s} \le \theta \le \frac{2\pi(k+1)}{s}}\\
0,\ \mbox{otherwise\,.}
\end{cases}
\end{equation}
Therefore, the number of boundary weight functions, denoted $s$ in Section \ref{s5}, will be $2s$ here.
The Lagrange multiplier has the form $\ds{\lambda^{(s)}=\sum_{k=0}^{s-1} \left(c_{k}\,g_k(\theta)+d_k\,h_k(\theta)\right)}$, and the Lagrange functional is $\ds{\mathscr{L}\Big(u^{(N)}\Big)=J(u^{(N)})+\int_{\partial\Omega}\lambda^{(s)}u^{(N)}\,d\sigma}$. The unknown coefficients $a_{i,j}$, $c_k$ and $d_k$ are now determined from the system of $N^2$ internal 
$\ds{\frac{\partial\mathscr{L}}{\partial a_{i,j}}=0}$ and $2s$ boundary $\ds{\int_{\partial\Omega}g_ku^{(N)}\,d\sigma=\int_{\partial\Omega}h_ku^{(N)}\,d\sigma=0}$ equations.

The relative errors of the Ritz-Lagrange solutions versus the exact solution $\ov{u}$ \eqref{CiSol} are shown in  Table 
\ref{cap:membrane_errors} as the percentages of the maximum deflection at $x=0$ and $y=0$ . They are quite small considering that one has to determine $N^2+2s$ coefficients in each case. Note that we always keep $N^2>2s$ as recommended in the description of the method to ensure that the system matrices have full rank and are invertible.
\renewcommand{\tablename}{\textbf{Table}}
\begin{table}[H]
\centering
\begin{tabular}{|c|c|c|c|}
\hline 
$$ & $\text{Central Error \%}$ & $\text{Boundary Error \%}$\tabularnewline
\hline 
$N=3, s=2$ & -4.04 & -1.33\tabularnewline
\hline 
$N=4, s=3$ & -4.05 & -0.1\tabularnewline
\hline 
$N=5, s=4$ & -0.04 & 0.03\tabularnewline
\hline 
\end{tabular}
\caption{\label{cap:membrane_errors} Central ($x=0, y=0$) and boundary ($x=0, y=1$) error relative to the maximum deflection of the membrane of unit radius.}
\end{table}
 
\textbf{Problem 5.} Consider a problem of bending a uniformly loaded, simply supported on all sides (SS-SS-SS-SS), isotropic, square plate of constant thickness, unit stiffness and unit edge length. Simply supported means that 
$u=0$ on $\partial{\Omega}$. We do not need to list the natural boundary conditions since a variational formulation incorporates them automatically. The variational functional giving the potential energy of the plate 
is \cite[IV.10.3]{CH}:
\begin{multline}\label{pot_energy}
J(u)=\int_{0}^{1}\int_{0}^{1}\frac{1}{2}\left(\left(\frac{\partial^{2}u}{\partial x^{2}}\right)^{2}+\left(\frac{\partial^{2}u}{\partial y^{2}}\right)^{2}+2\upsilon\,\frac{\partial^{2}u}{\partial x^{2}}\frac{\partial^{2}u}{\partial y^{2}}+2(1-\upsilon)\left(\frac{\partial^{2}u}{\partial x\partial y}\right)\right)-fu\,dxdy\,,
\end{multline}
where $u$ is the displacement of the plate, $\upsilon$ is the Poisson ratio, and $f$ is a distributed load. 

The Euler-Lagrange equation induced by Eq.\eqref{pot_energy} is the biharmonic equation $\nabla^2\nabla^2u=f$, the terms multiplied by the Poisson ratio form a divergence and only affect the natural boundary conditions.
As the internal trial functions we choose the products of cosines $X_i(x)=\cos((i-1)\pi x)$ and $Y_i(y)=\cos((i-1)\pi y)$, so that the trial solution is $\ds{u^{(N)}=\sum_{i=1}^N \sum_{j=1}^Nc_{ij}\,X_{i-1}(x)Y_{j-1}(y)}$. Obviously, 
the trial functions do not satisfy the boundary condition. The Lagrange functional is
\begin{multline}\label{prob5LagrFunct}
\mathscr{L}\Big(u^{(N)}\Big)=J\Big(u^{(N)}\Big)+\int_0^1 \lambda_1^{(s)}(x) u(x,0)\,dx\\+\int_0^1 \lambda_2^{(s)}(y) u(1,y)\,dy-\int_0^1 \lambda_3^{(s)}(x) u(x,1)\,dx-\int_0^1 \lambda_4^{(s)}(y) u(0,y)\,dy\,,
\end{multline}
where for convenience we split the Lagrange multiplier $\lambda^{(s)}$ into its restrictions $\lambda_i^{(s)}$ to each edge of the plate. This way we can represent the set of the boundary weight functions as the union of four sets selected separately for each edge, namely $\ds{\lambda_i^{(s)}(z)=\sum_{j=1}^s \lambda_{i,j} \cos\big((j-1)\pi z\big)}$, where 
$\lambda _{i,j}$ are the unknown coefficients. The number of internal equations here is again $N^2$, and the number of the boundary equations is $4s$, so the non-degeneracy condition is $N^2>4s$. 

The exact solution to this problem can be expressed as a rapidly convergent series, we use its first ten terms to calculate the errors. We do not tabulate them here, because they are very large (up to 70$\%$), and increasing the number of terms does not improve the approximation. At this point, the reader should not be surprised, indeed we are dealing with the same mistake as in the one-dimensional Problem 2. As we pointed out in Section \ref{s3}, the system of cosines is incomplete on the interval, so the system of their products naturally is incomplete on the product of intervals that represents the plate. A more down to Earth explanation is that products of cosines have vanishing normal derivatives 
$\ds{\,\frac{\partial u}{\partial n}}\,$ on all edges and all their linear combinations inherit this property. This does not matter for second order equations because the normal derivatives are discontinuous on the relevant spaces, but it does matter for the higher order equations like the biharmonic equation. 

One can check that in the weak formulation of Problem 5 because of the vanishing normal derivatives the boundary terms that multiply the variation of the solution's derivatives get removed, so we ended up solving a different variational problem. Indeed, the choice of cosine products unwittingly enforces an additional boundary condition, 
$\ds{\frac{\partial u}{\partial n}}=0$ on $\partial{\Omega}$. Together with $u=0$ on $\partial{\Omega}$ this describes, physically, a plate clamped on all sides (C-C-C-C) rather than a simply supported one. Thus, we should be comparing our Ritz-Lagrange solutions to the answers for the C-C-C-C plate (cf. Problem 3 in the one-dimensional case). Unfortunately, an analytic solution for a plate clamped on all sides is not known, so we used the values obtained in \cite[VI.44]{Tim} to make the comparison. The relative errors as percentages of the maximum deflection at the center of the plate are shown in Table \ref{cap:cccc_errors}. 
\begin{table}[H]
\centering
\begin{tabular}{|c|c|c|c|}
\hline 
$$ & $\text{Central Error \%}$ & $\text{Boundary Error \%}$\tabularnewline
\hline 
$N=4, s=2$ & 52.76 & -50.92\tabularnewline
\hline 
$N=6, s=3$ & -3.3 & -1.22\tabularnewline
\hline 
$N=8, s=4$ & 0.063 & -1.39\tabularnewline
\hline  
$N=10, s=5$ & -0.12 & -0.09\tabularnewline
\hline 
\end{tabular}
\caption{\label{cap:cccc_errors} Central ($x=0.5, y=0.5$) and boundary ($x=0, y=0$) error   relative to maximum bending deflection of an C-C-C-C square plate of unit size.}
\end{table}
 
Section \ref{s3} also gives us a way to solve the original problem, we just need to complete the system of cosine products. By Lemma \ref{MultDComp} it suffices to complete the cosines on the interval and take the products from the completed system. Namely, we take the products of $x,x^2,\cos((i-1)\pi x)$ and $y,y^2,\cos((i-1)\pi y)$ as the new trial functions, and keep the rest of the above setup intact. The relative errors for the Ritz-Lagrange solutions with the completed system against the known series solution \cite[8.2.4]{Reddy} are shown in Table \ref{cap:ssss_errors} and demonstrate the validity of the method.
\begin{table}[H]
\centering
\begin{tabular}{|c|c|c|c|}
\hline 
$$ & $\text{Central Error \%}$ & $\text{Boundary Error \%}$\tabularnewline
\hline 
$N=5, s=2$ & 26.14 & -89.03\tabularnewline
\hline 
$N=6, s=3$ & 2.81 & -3.7\tabularnewline
\hline 
$N=8, s=4$ & 1.68 & -4.4\tabularnewline
\hline 
$N=10, s=5$ & 0.56 & -1.1\tabularnewline
\hline 
\end{tabular}
\caption{\label{cap:ssss_errors} Central ($x=0.5, y=0.5$) and boundary ($x=0, y=0$) error   relative to maximum bending deflection of an SS-SS-SS-SS square plate of unit size.}
\end{table}

As a final demonstration, we apply the Ritz-Lagrange method to a boundary eigenvalue problem for square plates.  The eigenmodes describe standing vibrations of a plate, and their zeros (nodal curves) are known as Chladni figures 
\cite[5.1]{Gander}.  The problem has attracted a lot of attention from both analytic and numerical viewpoints, indeed Ritz himself applied his method to it in his original paper. Boundary eigenvalue problems are somewhat beyond the scope of the theory in Section \ref{s5}, which deals with linear constraints only, because under the Rayleigh-Ritz approach one needs to impose an additional quadratic normalization constraint $\frac12\int_\Omega u^2\,dx=1$ to solve for the eigenmodes \cite[18.5]{Zeid2a}, \cite[VI.1.1]{CH}, \cite[5.2]{Gander}. However, the general approach of Section \ref{s5} remains valid, and with some extra effort one can justify applying the Ritz-Lagrange method to problems with non-linear constraints along the same lines.

\textbf{Problem 6.}  Consider a uniformly loaded, simply supported on all sides, isotropic square plate of constant thickness with unit edge length. The potential energy $J$ of the plate is given by Eq.\eqref{pot_energy} without the distributed load term. The boundary eigenvalue problem can be interpreted as finding extrema of $J(u)$ subject to the boundary condition $u=0$ on $\partial{\Omega}$, and the normalization constraint $\frac12\int_\Omega u^2\,dx=1$.

Compared to Problem 5 the Lagrange functional acquires an additional term $\mu(\int_\Omega u^2\,dx-1)$ and an additional equation, which amounts to the normalization constraint on the eigenvectors. Of course, in practice one can ignore this equation and simply use standard methods for finding eigenvectors. We keep the choices for the internal and the boundary weight functions from Problem 5. Let $c$ denote the vector of internal coefficients $c_{ij}$ and $\l$ denote the vector of boundary coefficients $\l_{i,j}$. In terms of $c$ and $\l$ the Lagrange functional can be conveniently represented as $\mathscr{L}(u^{(N)})=\frac{1}{2}c^TKc-\mu\left(\frac{1}{2}c^TMc-1\right)+\l^TLc$, where $K$ and $M$ are matrices of size $N^2\times N^2$ obtained by integrating the internal trial functions, see \cite[8.2.7]{Reddy}, and $L$ is a $4s\times N^2$ matrix obtained by integrating the boundary weight functions. Matrix $L$ can be obtained by multiplying the boundary equations with the corresponding Lagrange multiplier functions and extracting the coefficients of  $c_{ij}$ and 
$\l_{i,j}$ after the integration. We note that the boundary equations can be written as $Lc=0$. Finally, differentiating the Lagrangian with respect to $c_{ij}$ and $\l_{i,j}$ we are led to the following generalized eigenvalue problem:
\begin{equation}
\left(\left[\begin{array}{c|c}
K\ & \ L^{T}\\ \hline
L\ & \ 0
\end{array}\right]-\mu\left[\begin{array}{c|c}
M\ & \ 0\ \\ \hline
0\ & 0
\end{array}\right]\right)
\left(\begin{array}{c}c\\\lambda\end{array}\right)=0
\end{equation}
For this eigenvalue problem to be solvable one needs $L$ to have the maximal rank $4s$, which is ensured by the non-degeneracy condition $N^2\gg 4s$ as in Problem 5. The eigenvalues $\mu_i=\omega_i^2$ approximate the squares of the natural frequencies of the plate's vibrations.

With $N=10$ and $s=5$ we obtain a set of approximate non-dimensional natural frequencies $\omega_i$, first nine of which are shown below. Since the eigenmodes are known to be of the form $\sin(\pi m x)\sin(\pi n y)$ we change the single index notation to $\omega_{mn}$ and arrange the frequencies in a square pattern
\begin{equation}
\left[\begin{array}{ccc}
19.61797 & 49.06479 & 98.33527\\
49.06479 & 77.55724 & 126.62091\\
99.03778 & 126.62091 & 177.73211
\end{array}\right]\,.
\end{equation}
The exact values are taken from \cite[8.2.4]{Reddy}, repeated frequencies corrrespond to multiple eigenvalues with the eigenmodes symmetric along different axes:
\begin{equation}
\left[ \begin {array}{ccc}  19.73920881& 49.34802202& 98.69604404
\\ \noalign{\medskip} 49.34802202& 78.95683523& 128.3048573
\\ \noalign{\medskip} 98.69604404& 128.3048573& 177.6528793
\end {array} \right]\,.
\end{equation}
One can see that the estimated frequencies are slightly lower than the exact ones. This is in contrast with the application of the usual Ritz method, where the estimated frequencies are always higher. From a physical viewpoint, the latter happens because replacing an infinite system with a finite one is equivalent to imposing additional constraints, which tend to raise the stiffness of the system, and hence the frequencies. This assumes however that all the boundary constraints are enforced in both systems, i.e. that the trial functions satisfy the essential boundary conditions. 

In the Ritz-Lagrange method the trial functions do not satisfy the essential conditions, and even 
the trial solutions are forced to satisfy them only approximately. In other words, i.e. {\it we are effectively relaxing the boundary constraints in addition to imposing additional ones through discretization. This relaxation lowers the frequencies (because a plate with fewer constraints is less stiff) and counteracts the effects of discretization}. If we were able to impose the boundary conditions everywhere along the boundary the estimated frequencies would have been higher than the exact ones just as in the usual Ritz method. 

From a mathematical viewpoint, this effect is also natural since the eigenvalues are the minima of a quadratic functional on subspaces of the original space \cite[VI.1.1]{CH}. In the Ritz-Lagrange method we approximate them by using functions from a larger space (by relaxing the boundary conditions), thus lowering the minima that can be attained. In particular, one can see from the proof of Theorem \ref{RitzLagBnd} that the values of the functional at the approximating elements are potentially smaller than at the sought minimizer.

\section{Conclusions}\label{s7}

We developed a general extension of the Ritz method to systems of trial functions that do not satisfy the essential boundary conditions, and proved its convergence. The method is based on treating the essential conditions as variational constraints and removing them using the Lagrange multipliers. Here are some general observations of the workings of the method.
\begin{itemize}
\item The variational functional has to be well-behaved not only on the energy space of the problem, but on its extension that contains the trial functions. Sufficiently good behavior is a strong form of convexity, which in the case of quadratic functionals means that the boundary value problem is strongly elliptic.
\item The systems of trial functions must be complete in the norms consistent with the functional, which usually restrict to the energy norms on the energy space of the problem. Although similar requirement applies to the usual Ritz method, here it is much easier to encounter systems that appear complete but are not due to effects at the boundary.
\item The Lagrange multipliers have to be treated as additional variables in the approximating systems. They can not be eliminated by substituting the trial solutions into the variational formulas for them in terms of the exact solution. These formulas are discontinuous in the relevant norms.
\item In multidimensional problems the boundary conditions incorporate infinitely many constraints, and to obtain a finite dimensional approximating system one has to select boundary weight functions in addition to the trial functions. The number of trial functions has to be significantly larger than the number of the boundary weights, otherwise the approximating system may be inconsistent or only have the trivial solution.
\item In multidimensional problems the approximating values of the functional may approach the exact value from below rather than above, as in the usual Ritz method, because the minimization takes place on a larger space of functions not satisfying the boundary conditions exactly.
\item The method can be applied to boundary eigenvalue problems interpreted along the lines of Rayleigh-Ritz as minimization problems on subspaces of the original space with the additional normalization constraint. Due to the presence of the Lagrange multiplier variables the resulting finite dimensional problem is a generalized eigenvalue problem $(A-\mu B)x=0$ instead of the ordinary one with $B=I$. In multidimensional vibrational problems the approximate eigenfrequencies obtained in this way may be lower than the exact ones, in contrast to the Ritz method where they are always higher, due to relaxation of the boundary constraints.
\end{itemize}

As is well-known \cite{Leip,Reddy}, the Ritz method leads to the same approximating systems as the Galerkin method, but the latter can also be applied to non-optimization problems. It is quite intriguing whether the Ritz-Lagrange method developed here can be extended to a 'Galerkin-Lagrange method' for non-variational problems. There is no Lagrange functional to be had in such problems, but one can formally add boundary terms multiplied by extra variables ('Lagrange multipliers') to a weighted residual of the problem. An approximation of the boundary conditions should also be added to the usual Galerkin system. However, as we saw in the Lanczos tau example such a straightforward approach is not likely to work. Indeed, in the Ritz-Lagrange method we add the Lagrange boundary terms not to the bare weighted residual, but to an integrated by parts expression with some boundary terms of its own. This suggests that in a correct generalization the problem has to be rewritten in a weak form \cite[7.5.1]{Reddy} before the Lagrange-like boundary terms are added.

Another complication is the role of the natural boundary conditions. In variational problems they are enforced automatically, so there is no point in adding them as constraints and introducing additional Lagrange multipliers. An example in \cite{Par} even shows that attempting to do so leads to worsening the convergence of the trial solutions. However, it is unclear how the natural conditions should be enforced in a 'Galerkin-Lagrange method'. Still, the distinction between the natural and the essential boundary conditions in \cite[I.1.2]{Col} (by the order of the derivatives entering them) makes sense even for non-optimization problems, and it has been shown in \cite{Leip} that in some cases the weak form of a problem provides an enforcement mechanism for the natural conditions.

\section*{Appendix: Proofs}
\setcounter{equation}{0} 
\renewcommand{\theequation}{A.\arabic{equation}}

\begin{proof}[Proof of {\bf Lemma \ref{U0comb} }] 
Let $\{\phi_i\}$ be a complete system in $\U$. Since $\G$ has an $s$-dimensional image we can represent it as $\G u=\big(\langle \G_1,u\rangle,\dots,\langle\G_s,u\rangle\big)^T$, where $\G_i$ are bounded linear functionals. Assume without loss of generality that they are linearly independent, otherwise some of them can be dropped without changing $\o{\U}$. Set $\o{\U}_k:=\{u\in\U\,\big|\,\langle \G_1,u\rangle,\dots,\langle \G_k,u\rangle=0\}$, we will construct a complete system in each $\,\o{\U}_k$ by induction on $k$. Since $\,\o{\U}=\o{\U}_s$ the process concludes in $s$ steps.

For $k=1$ we must produce a complete system of linear combinations in 
$\o{\U}_1:=\{u\in\U\,\big|\,\langle \G_1,u\rangle=0\}$. Without loss of generality, 
$\langle \G_1,\phi_1\rangle\neq0$ since $\{\phi_i\}$ is complete and $\G_1$ can not vanish on all $\phi_i$. 
We claim that 
$\ds{\widetilde{\phi}_i:=\phi_i-\frac{\langle \G_1,\phi_i\rangle}{\langle \G_1,\phi_1\rangle}\,\phi_1\in\o{\U}_1}$ 
form the desired system. Let $u\in\o{\U}_1\subset\U$ and $a_i$ be the coefficients such 
that $\|u-\sum_{i=1}^{N}a_i\,\phi_i\|\leq\varepsilon$ for a given $\varepsilon>0$. By definition of 
$\widetilde{\phi}_i$,
$$
\sum_{i=1}^{N}a_i\,\widetilde{\phi}_i=\sum_{i=1}^{N}a_i\,\phi_i
-\frac{\langle \G_1,\sum_{i=1}^{N}a_i\,\phi_i\rangle}{\langle \G_1,\phi_1\rangle}\,\phi_1\,.
$$
To estimate the second term we find,
$$
|\langle\G_1,\sum_{i=1}^{N}a_i\,\phi_i\rangle|=
|\langle\G_1,\sum_{i=1}^{N}a_i\,\phi_i-u\rangle+\langle \G_1,u\rangle|
\leq\|\G_1\|\,\|u-\sum_{i=1}^{N}a_i\,\phi_i\|\leq\|\G_1\|\varepsilon\,.
$$
Therefore, $\ds{\|u-\sum_{i=1}^{N}a_i\,\widetilde{\phi}_i\|
\leq\Big(1+\frac{\|\G_1\|\,\|\phi_1\|}{|\langle \G_1,\phi_1\rangle|}\Big)\varepsilon}$, and since $u$, $\varepsilon$ are arbitrary completeness of $\widetilde{\phi}_i$ follows.

Let $\{\widetilde{\phi}_i\}$ be a complete system in $\o{\U}_k$ from the preceeding step. 
Linear independence of $\G_j$ guarantees that $\G_{k+1}$ does not vanish on some $\widetilde{\phi}_i$, which we may as well take to be $\widetilde{\phi}_1$. Apply the process above with 
$\G_1$ replaced by $\G_{k+1}$ and $\o{\U}_1$ replaced by $\o{\U}_{k+1}$ to obtain $\widetilde{\widetilde{\phi}}_i$. 
Then $\widetilde{\widetilde{\phi}}_i$ are linear combinations of $\widetilde{\phi}_i$ 
(and hence of the original $\phi_i$ ), belong to $\o{\U}_{k+1}$ and are complete in it by the same argument. 
This concludes the induction step.
\end{proof}

\begin{proof}[Proof of {\bf Theorem \ref{RitzLagCon}}]
A standard argument from convex analysis shows that if $J(u)\xrightarrow[\|u\|\to\infty]{}\infty$ on $\o{\U}$ then $J$ has minimizers on $\o{\U}$ and there is a weakly convergent subsequence $u^{\,(N_k)}\xrightarrow[k\to\infty]{w}u^{\,(\infty)}$ \cite[III.10.3]{Vain}, \cite[6.2]{DraMil}. For convex functionals strong continuity implies weak continuity so $J\big(u^{\,(N_k)}\big)\xrightarrow[k\to\infty]{}J\big(u^{\,(\infty)}\big)$, and moreover $u^{\,(\infty)}\in\o{\U}$ since $0=\G\big(u^{\,(N_k)}\big)\xrightarrow[k\to\infty]{}\G\big(u^{\,(\infty)}\big)$. But for large enough $N$ there is a $u\in\o{\U}^{\,(N)}$ arbitrarily close to a minimizer $\ov{u}$ of $J$ on $\o{\U}$ by Lemma \ref{U0comb}, hence by continuity $J(u)$ is arbitrarily close to the minimal value $J_{\min}$. But $J\big(u^{\,(N_k)}\big)$ can not exceed $J(u)$ for $N_k\geq N$ since it is a minimizer on $\o{\U}^{\,(N_k)}$, so $J_{\min}\leq J\big(u^{\,(N_k)}\big)\leq J(u)=J_{\min}+\varepsilon$. After passing to limit we have that $J\big(u^{\,(\infty)}\big)=J_{\min}$, i.e. $u^{\,(\infty)}$ is a minimizer of $J$ on $\o{\U}$. 

If we assume additionally that $J$ is strictly convex on $\o{\U}$ then $\ov{u}$ is unique and the entire sequence $u^{\,(N)}$ (which is now also uniquely defined) converges to it at least weakly \cite[6.2A]{DraMil}.
\end{proof}

The next two proofs use equivalent norms (inner products) on $W_2^1([0,\pi])$ and $W_2^2([0,\pi])$ respectively. Two norms are equivalent if they define the same notion of convergence, for equivalent norms on Sobolev spaces see 
\cite[I.8]{Lad} and especially \cite[1.9]{Sob}.
\begin{proof}[Proof of {\bf Lemma \ref{CosW21}}] 
The following inner product is equivalent to the usual one on $W_2^1([0,\pi])$:
$\ds{\langle u,v\rangle_0:=u(0)v(0)+\int_0^\pi u_xv_x\,dx}$. To prove completeness it suffices to show that any 
function $w$ orthogonal to all cosines must be $0$. For such $w$ we have $\langle w,1\rangle_0=w(0)=0$ and 
hence $\ds{\langle w,\cos nx\rangle_0=\int_0^\pi w_x\cdot(-n\sin nx)\,dx=0}$ for $n\geq1$. Thus, $w_x$ is $L_2$ orthogonal to $\sin nx$ for all $n\geq1$. Since the latter form an orthogonal basis in $L_2([0,\pi])$ we must have 
$w_x=0$ a.e. But then by the Fundamental Theorem of Calculus $\ds{w(x)=w(0)+\int_0^x w_t\,dt=0}$ a.e. establishing completeness. Being an orthogonal basis in $L_2([0,\pi])$ cosines must be minimal there, and 
therefore in any space with a stronger norm, which includes $W_2^1([0,\pi])$.
\end{proof}

\begin{proof}[Proof of {\bf Lemma \ref{CosW22}}] An equivalent inner product on $W_2^2([0,\pi])$ is\\
$\ds{\langle u,v\rangle_0:=u(0)v(0)+u_x(0)v_x(0)+\int_0^\pi u_{xx}v_{xx}\,dx}$. Consider $w$ orthogonal to all cosines, then we have $\langle w,1\rangle_0=w(0)=0$ and $\ds{\langle w,\cos nx\rangle_0=\int_0^\pi w_{xx}\cdot(-n^2\cos nx)\,dx=0}$ for $n\geq1$ because all sines vanish at $0$. In particular, $w_{xx}$ is $L_2$ orthogonal to $\cos nx$ for all $n\geq1$. But orthogonal complement of the latter in $L_2$ consists of constants, so $w_{xx}=\text{const}$ and $w(x)=ax^2+bx+c$. Since 
$w(0)=0$ free term is $0$ and $w$ is a linear combination of $x$ and $x^2$. Thus, orthogonal complement to cosines is spanned by $x$ and $x^2$ proving completeness.

For minimality notice that by direct calculation 
$\langle x,\cos nx\rangle_0=\langle x^2,\cos nx\rangle_0=\langle x,x^2\rangle_0=0$, i.e. $x$ and $x^2$ are orthogonal to all cosines and to each other. This means that neither one of them can be deleted without loosing completeness. It also means that if a cosine can be approximated in $W_2^2$ by other cosines combined with $x$ and $x^2$ then it can already be approximated by other cosines alone. But the latter can not be done with arbitrary precision even in $L_2$, let alone in $W_2^2$.
\end{proof}

\begin{proof}[Proof of {\bf Theorem \ref{RitzLagBnd}}]
Since $\o{\U}\subset\cdots\subset\o{\U}_2\subset\o{\U}_1$ and the minimum on a larger space can not get bigger we have $J(\ov{u})\geq\cdots J(\ov{u}_2)\geq J(\ov{u}_1)$. Thus, the numerical sequence $J(\ov{u}_s)$ is bounded. Moreover, $\ov{u}_s\in\o{\U}_{s_0}$ for $s\geq s_0$, so 
$\|\ov{u}_s\|\leq M<\infty$ for $s\geq s_0$ since $J$ grows at infinity on $\o{\U}_{s_0}$. Recall that $\ov{u},\,\ov{u}_s$are the minimizers of  $J$ on $\o{\U},\,\o{\U}_s$ respectively, and therefore the derivatives $J'(\ov{u}),\,J'(\ov{u}_s)$ vanish when paired with elements of the corresponding subspaces. In particular, 
$\langle J'(\ov{u}),\ov{u}\rangle=\langle J'(\ov{u}_s),\ov{u}\rangle=\langle J'(\ov{u}_s),\ov{u}_s\rangle=0$ and 
\begin{equation}\label{Umonus}
\langle J'(\ov{u}_s)-J'(\ov{u}),\ov{u}_s-\ov{u}\rangle=-\langle J'(\ov{u}),\ov{u}_s\rangle.
\end{equation}
We will prove that the last expression converges to $0$ when $s\to\infty$. Condition Eq.\eqref{Umon} then implies
$c\,(\ov{u}_s-\ov{u})\xrightarrow[s\to\infty]{}0$ and hence $\ov{u}_s\xrightarrow[s\to\infty]{}\ov{u}$ as claimed. Convergence of $J(\ov{u}_s)$ follows from continuity of $J$.

We now prove convergence in  Eq.\eqref{Umonus}. Since $\ov{u}$ is a minimizer on $\o{\U}$ the functional $J'(\ov{u})$ vanishes on any element from it. The subspace of functionals that vanish on the entire 
$\o{\U}=\{u\in\U\,\big|\,\G u=0\}$ is the closed linear span of $\{\G^*\psi_j\}$ in $\U^*$. Indeed,  if $\{\G^*\psi_j\}$ did not span it there would exist, by the Khan-Banach theorem, a $u$ such that $\langle \G^*\psi_j,u\rangle=\langle\psi_j,\G u\rangle=0$ for all $j$, while $\G u\neq0$, contradicting the completeness of $\{\psi_j\}$.
Thus, for any $\varepsilon>0$ there exists a linear combination $\xi=\sum_{j=1}^na_j\G^*\psi_j$ such that $\|J'(\ov{u})-\xi\|\leq\varepsilon$. But then $\xi\in\o{\U}_n$, and for $s>n$ we have $\langle\xi,\ov{u}_s\rangle=0$, so\\ 
\centerline{$|\langle J'(\ov{u}),\ov{u}_s\rangle|\leq\|J'(\ov{u})-\xi\|\,\|\ov{u}_s\|\leq M\varepsilon$.}
Since $\varepsilon$ is arbitrary $\langle J'(\ov{u}),\ov{u}_s\rangle\xrightarrow[s\to\infty]{}0$.
\end{proof}

\begin{proof}[Proof of {\bf Lemma \ref{MultDComp}}] In the multiindex notation an equivalent norm on $W_p^k(\D)$ is
given by
$$
\|F\|_{W_p^k}:=\sum_{i=0}^k\sum_{|\alpha|=i}\left\|\frac{\d^{|\alpha|}F}{\d\xi^\alpha}\right\|_{L_p}
=\sum_{|\alpha|\leq k}\left\|\frac{\d^{|\alpha|}F}{\d\xi^\alpha}\right\|_{L_p}\!\!\!\!, 
$$
where the $L_p$ norm is just $\|f\|_{L_p}:=(\int_\D|f|^p\,d\xi)^{\frac1p}$. If $f\in W_p^k(\O)$ and 
$\widetilde{f}\in W_p^k(\widetilde{\O})$ then it follows from the Fubini theorem that 
$\|f\widetilde{f}\|_{L_p}=\|f\|_{L_p}\|\widetilde{f}\|_{L_p}$ since $f$ and $\widetilde{f}$ depend 
on different variables. Let $x$ and $\widetilde{x}$ denote the variables on $\O$ and $\widetilde{\O}$ respectively, so that $\xi=(x,\widetilde{x})$ is the variable on $\O\times\widetilde{\O}$. Then we estimate
\begin{multline*}
\|F\widetilde{F}\|_{W_p^k}
=\sum_{|\alpha|\leq k}\Big\|\frac{\d^{|\alpha|}F\widetilde{F}}{\d\xi^\alpha}\Big\|_{L_p}
=\sum_{|\beta|+|\gamma|\leq k}\Big\|\frac{\d^{|\beta|}F}{\d x^\beta}\Big\|_{L_p}
\Big\|\frac{\d^{|\gamma|}\widetilde{F}}{\d\widetilde{x}^\gamma}\Big\|_{L_p}\\
\leq\sum_{|\beta|\leq k,|\gamma|\leq k}\Big\|\frac{\d^{|\beta|}F}{\d x^\beta}\Big\|_{L_p}
\Big\|\frac{\d^{|\gamma|}\widetilde{F}}{\d\widetilde{x}^\gamma}\Big\|_{L_p}
=\sum_{|\beta|\leq k}\Big\|\frac{\d^{|\beta|}F}{\d x^\beta}\Big\|_{L_p}
\sum_{|\gamma|\leq k}\Big\|\frac{\d^{|\gamma|}\widetilde{F}}{\d\widetilde{x}^\gamma}\Big\|_{L_p}
=\|F\|_{W_p^k}\|\widetilde{F}\|_{W_p^k}\,.
\end{multline*}
Since $\phi_i$ are complete any monomial $x^\beta$ can be approximated to any precision $\varepsilon>0$ in $W_p^k$ by their linear combination $\phi=\sum_{i}a_i\,\phi_i$, and analogously $\widetilde{x}^\gamma$ can be approximated by a linear combination $\widetilde{\phi}=\sum_{j}\widetilde{a}_j\,\widetilde{\phi}_j$. But $\phi\widetilde{\phi}=\sum_{i,j}a_i\widetilde{a}_j\,\phi_i\widetilde{\phi}_j$ is a linear combination of $\phi_i\widetilde{\phi}_j$, while the difference 
between the products can be made arbitrarily small:
\begin{multline*}
\|x^\beta\widetilde{x}^\gamma-\phi\widetilde{\phi}\|_{W_p^k}
=\|x^\beta(\widetilde{x}^\gamma-\widetilde{\phi})+(x^\beta-\phi)\widetilde{\phi}\|_{W_p^k}\\
\leq\|x^\beta\|_{W_p^k}\|\widetilde{x}^\gamma-\widetilde{\phi}\|_{W_p^k}
+\|x^\beta-\phi\|_{W_p^k}\|\widetilde{\phi}\|_{W_p^k}
\leq\varepsilon\,(\|x^\beta\|_{W_p^k}+\|\widetilde{x}^\gamma\|_{W_p^k}+\varepsilon)\,, 
\end{multline*}
\noindent where the first inequality follows from the above estimate. Hence any product of monomials, and therefore any polynomial, can be approximated in $W_p^k$ by linear combinations of 
$\phi_i\widetilde{\phi}_j$. By the generalized Weierstrass theorem \cite[II.4.3]{CH}, polynomials are complete in 
$W_p^k(\O\times\widetilde{\O})$, and hence so is the system $\{\phi_i\widetilde{\phi}_j\}$.
\end{proof}

\bibliographystyle{unsrt}
\bibliography{RitzLag}

\end{document}